\documentclass[lettersize,journal]{IEEEtran}
\usepackage{amsmath,amsfonts}
\usepackage{algorithmic}
\usepackage{algorithm}
\usepackage{array}
\usepackage[caption=false,font=normalsize,labelfont=sf,textfont=sf]{subfig}
\usepackage{textcomp}
\usepackage{stfloats}
\usepackage{url}
\usepackage{verbatim}
\usepackage{graphicx}
\usepackage{cite} 
\usepackage{amssymb}
\usepackage{color}
\usepackage{enumitem}
\usepackage{multirow}
\usepackage{booktabs}
\allowdisplaybreaks

\begin{document}

\title{A Distributionally Robust Optimization Framework for Stochastic Assessment of Power System Flexibility in Economic Dispatch}

\author{Xinyi Zhao, Lei Fan, Fei Ding, Weijia Liu, and Chaoyue Zhao
\thanks{Xinyi Zhao and Chaoyue Zhao  (\{xyzhao24; cyzhao\}@uw.edu) are with the Department of Industrial \& Systems Engineering, University of Washington, Seattle WA, USA. Lei Fan (lfan8@central.uh.edu) is with the Department of Engineering Technology, University of Houston, Houston TX, USA. Fei Ding and Weijia Liu  (\{Fei.Ding; Weijia.Liu\}@nrel.gov) are with the National Renewable Energy Laboratory, Golden CO, USA.}
\thanks{This work is funded by NSF ECCS 2046243, 2045978.}}


\maketitle

\begin{abstract}
Given the complexity of power systems, particularly the high-dimensional variability of net loads, accurately depicting the entire operational range of net loads poses a challenge. To address this, recent methodologies have sought to gauge the maximum range of net load uncertainty across all buses. In this paper, we consider the stochastic nature of the net load and introduce a distributionally robust optimization framework that assesses system flexibility stochastically, accommodating a minimal extent of system violations. We verify the proposed method by solving the flexibility of the real-time economic dispatch problem on four IEEE standard test systems. Compared to traditional deterministic flexibility evaluations, our approach consistently yields less conservative flexibility outcomes.
\end{abstract}


\section{Introduction}
The worldwide electricity sector has witnessed a significant rise in renewable energy integration. As this trend is projected to continue in the coming decades, the net demand, calculated by subtracting electricity generation from the total load, will become more volatile and unpredictable. 
To address this issue cost-effectively, system operators need to evaluate the flexibility of the power system, which describes the system's capability of managing the variability and uncertainty of net loads \cite{lannoye2012evaluation, ma2013evaluating}. It is critical to gain insights into flexibility as it aids operators in anticipating unexpected shifts in demand. Moreover, overlooking the importance of flexibility could result in transient instabilities, cascading failures, and the risk of widespread blackouts.

Reviews \cite{lannoye2012evaluation,lannoye2012power,mohandes2019review} on power system flexibility categorize prior studies into two main groups based on the time scope of their target applications. First, from a short-term operational viewpoint, system frequency is a pivotal indicator of electrical power quality. Deviations from its nominal value typically result from imbalances between generation and load, necessitating timely compensation from available resources. As such, the capacity for regulation, power range, and ramping duration serves as their flexibility indices \cite{ma2013evaluating}. However, these studies have inherent limitations: a flexibility index relevant in one context might be inapplicable in another. In contrast, defining flexibility based on system failure causes, rather than remedies, makes it more universally relevant across different operational systems, offering greater utility for power grid operators.

Second, from a long-term planning perspective, scholars have introduced diverse technical and economic indices to gauge system flexibility from multiple facets. These include generation adequacy metrics such as the loss of load expectation \cite{calabrese1947generating}; ramping resource sufficiency indicated by the insufficient ramping resource expectation \cite{lannoye2014transmission}; and flexibility endurance, i.e., periods of flexibility deficit \cite{lannoye2012assessment}. These indices typically originate from simulations with preset net load probability assumptions. However, the task of computing the multi-dimensional joint probability distribution of net loads in large-scale, real-world power systems presents a formidable computational challenge. Furthermore, while these indices do well to capture system failures during certain variability patterns, such as Gaussian-distributed net loads, they fail to
adequately represent the operational range of net loads where systems function normally.

To address these challenges, state-of-the-art techniques concentrate on quantifying the utmost net load uncertainty that a system can accommodate \cite{shao2017security, zhao2015unified, wang2016robust}. These methodologies normally employ a two-stage robust model to gauge the worst-case scenarios arising from renewable energy unpredictability \cite{jiang2011robust, bertsimas2012adaptive}. Instead of optimizing for cost-related objectives, some strategies apply a robust optimization model to deterministically identify the maximum net load deviation from its typical baseline \cite{pulsipher2019computational, fan2021flexibility}. Nonetheless, such a conservative deviation range doesn't invariably ensure total operational safety, particularly given the rarity of worst-case scenarios. Typically, system operators are willing to tolerate minor levels of potential disruptions if they result in enhanced system adaptability. To cater to this perspective, our proposed stochastic assessment model evaluates system flexibility less conservatively, presenting two key contributions:
\begin{itemize}
\item The model identifies the maximum net load variation across all buses, ensuring that the expected operational violations stay within a predefined acceptable threshold for any net load profile within this variation range.

\item The model is employed in sequential real-time economic dispatch, independently assessing system flexibility at each time interval with hyperbox metrics. Our experiments demonstrate the advanced performance of our stochastic assessment in comparison to traditional deterministic approaches.
\end{itemize}

\section{Assessment Methods}

\subsection{Flexibility Metric}
The hyperbox metric evaluates the safe operating range, $\mathcal{U}$, of net loads in the power system \cite{pulsipher2019computational,fan2021flexibility}. Defined by $\Delta d= [\Delta d_b]$ as the peak deviation of the net load $d_b$ for each bus $b$ within the set $\mathcal{B}$, and $\bar{d}$ as the average or user-defined normal net loads. These parameters can be empirically derived from historical observations. The hyperbox representation of the uncertainty set is thus given as:
\begin{equation}
    U(\lambda) = \{\xi: \bar{d} - \lambda\Delta d \leq \xi \leq \bar{d} + \lambda\Delta d \}. \label{hyperbox}
\end{equation}
In \eqref{hyperbox}, a higher generic value $\lambda \in [0,1]$ indicates the system flexibility, while $\xi\in\mathbb{R}^{|\mathcal{B}|}$ is a realization of the random net load. This metric guarantees the minimum level of net load uncertainty tolerance. It is manifest that the flexibility set $U(\lambda) \subseteq \mathcal{U}$. A net load combination $\xi \notin U(\lambda) $ does not necessarily trigger a system failure. 

\subsection{Deterministic Assessment}\label{sec:det_ass}
The goal of the deterministic assessment for system flexibility is to identify the largest feasible $\lambda$ based on the flexibility metric such that the system can accommodate all $\xi \in U(\lambda)$. Let $\bold{x}$ be a vector including all decision variables. We then propose the following general optimization framework:
\begin{subequations}\label{linearabstract}
\begin{align}
\max \ \ & \lambda \\
\text{s.t.} \ \ &\max_{\xi\in U(\lambda)} \phi(\xi) \leq 0, \label{DRO_con} 
\end{align}
\end{subequations}
where
\begin{subequations}\label{phifunction}
\begin{align}
\phi(\xi)=\min_{\bold{x},u} \ \ & {\bf 1}^T u \\
\text{s.t.} \ \ & A_1\bold{x} - u_1\leq h_1 + H_1\xi,\label{RO_c1}\\
& A_2\bold{x} +u_2^+-u_2^- = h_2 + H_2\xi, \label{RO_c2}\\
& u \geq 0.
\end{align}
\end{subequations}

In the model, constraint \eqref{RO_c1} represents all system inequality constraints, whereas \eqref{RO_c2} captures all system equalities. The term $u$ denotes system violations, and the objective is to determine the maximum deviation $\lambda$ ensuring no system violations, even under the worst case $\xi$ running within $U(\lambda)$, as indicated by constraint \eqref{DRO_con}. A detailed mathematical model of this concept in \eqref{linearabstract} for the power system's economic dispatch problem will be introduced in Section \ref{sec:ED}.

\vspace{0.2cm}
\noindent {\bf Solution Approach:} Several methods have been proposed to address problem \eqref{linearabstract}. For example, \cite{pulsipher2019computational} establishes that maximizing $\lambda$ is equivalent to solving a mixed-integer program that reformulates the constraint $\phi(\xi) = 0$ using its first-order Karush-Kuhn-Tucker conditions. Furthermore, the cutting plane method in \cite{fan2021flexibility} also presents an alternative solution for addressing problem \eqref{linearabstract}.

\subsection{Stochastic Assessment}\label{sec:sto_ass}
Building upon the deterministic assessment, we extend our methodology \cite{fan2021flexibility} to develop a stochastic one, aimed at characterizing the uncertainty of the net load. We assume that $\xi$ follows a probability distribution denoted as $P(\xi)$, which belongs to the following ambiguity set:
\begin{equation}
  \mathcal{D}(\lambda) = \bigg\{P(\xi) \bigg| \int_{\xi \in U(\lambda)} dP(\xi) = 1, \ \int_{\xi \in U(\lambda)} \xi dP(\xi) = \bar{d}\bigg\}. \label{ambugity_set}
\end{equation}
This ambiguity set indicates that we consider all distribution $P(\xi)$ if its support is on $U(\lambda)$ and the mean value is $\bar{d}$. 

The primary objective of the stochastic assessment for system flexibility is to identify the most extensive support set within the ambiguity set $\mathcal{D}(\lambda)$, ensuring that the expected constraint violation, considering the worst-case distribution within $\mathcal{D}(\lambda)$, remains below a predefined threshold $\beta$. The abstract formulation can be expressed as follows:
\begin{subequations}\label{DRO_model}
\begin{align}
\max_{0\le\lambda\le 1} \ \ & \lambda \label{DRO_OF}\\
\text{s.t.} \ \ & \max_{P(\xi) \in \mathcal{D(\lambda)}}\ \  \mathbb{E}_{P(\xi)}[\phi(\xi)] \leq \beta. \label{DRO_CONS}
\end{align}
\end{subequations}

The resulting formulation \eqref{DRO_model} is a distributionally robust optimization (DRO) model \cite{delage2010distributionally}. In this variant, the distribution of the random parameter $\xi$ is uncertain and can vary adversely within the decision-dependent (endogenous) ambiguity set $\mathcal{D}(\lambda)$, with the optimal solution determined by considering the worst-case distribution.

To tackle \eqref{DRO_model}, we can treat the objective in \eqref{DRO_OF} as the master problem and redefine the internal maximization function within the constraints in \eqref{DRO_CONS} as the subproblem. Employing the ambiguity set specified in \eqref{ambugity_set}, we reformulate the maximization function in \eqref{DRO_CONS} as follows, where we represent $\mathbb{E}_{P(\xi)}[\phi(\xi)]$ as $\int_{\xi\in U(\lambda)} \phi(\xi) dP(\xi)$.

\begin{equation}\label{betaprob}
\begin{split}
\max_{P(\xi)} \ \ \bigg\{ &\int_{\xi\in U(\lambda)} \phi(\xi) dP(\xi): \int_{\xi\in U(\lambda)}dP(\xi)=1,\\
&\int_{\xi\in U(\lambda)}\xi dP(\xi) =\bar{d}. \bigg\}
\end{split}
\end{equation}
Let $\alpha$ and $\gamma$ serve as dual variables of two constraints in \eqref{betaprob}, its dual formulation can be expressed as follows:

\begin{equation}\label{dualsub}
\begin{split}
\min_{\alpha, \gamma} \ \ & \alpha+\bar{d}^{T}\gamma \\
\text{s.t.} \ \ & \alpha+\xi^{T}\gamma \geq \phi(\xi), \ \ \forall \xi\in U(\lambda),\\
& \alpha,\gamma\ \text{free}.
\end{split}
\end{equation}
Using the minimax duality for the Lagrangian, \eqref{dualsub} is equivalent to:

\begin{equation}\label{minmax}
\min_{\gamma} \ \ \left\{ \bar{d}^{T}\gamma \ +\max_{\xi\in U(\lambda)}\ \ (\phi(\xi) - \xi^{T}\gamma)  \right\}.
\end{equation}

To further express $\phi(\xi)$, we develop the dual formulation of the formulation \eqref{phifunction}:

\begin{subequations}\label{dualmas}
\begin{align}
\max_{\xi \in U(\lambda), \mu, \nu} \ \ & (h_1+H_1 \xi)^T \mu + (h_2+H_2 \xi)^T \nu \\
\text{s.t.} \ \ & A_1^T\mu + A_2^T\nu \le 0, \label{mu_c1}\\
& -{\bf 1}\le \mu\le {\bf 0}, \ -{\bf 1}\le \nu\le {\bf 1}. \label{mu_c2}
\end{align}
\end{subequations}
Here, $\mu$ and $\nu$ are introduced as dual variables for constraints \eqref{RO_c1} and \eqref{RO_c2}, respectively. Subsequently, we substitute $\phi(\xi)$ in \eqref{minmax} with the objective function derived in \eqref{dualmas}. This reformulation of the maximization function in \eqref{DRO_CONS} is presented as follows:

\begin{equation}\label{dualminmax}
\begin{split}
\min_{\gamma} \ \ \bigg[ &\bar{d}^{T}\gamma \ + \max_{\xi \in U(\lambda), \mu, \nu}\ \  \Big\{(h_1+H_1 \xi)^T \mu + (h_2+H_2 \xi)^T \nu \\
&- \xi^{T}\gamma: \text{Constraints} \ \eqref{mu_c1}-\eqref{mu_c2}\Big\}\bigg]\le \beta.
\end{split}
\end{equation}
Considering that ``min" in \eqref{dualminmax} indicates feasibility, it can be safely omitted in this context. As a result, the minimax formulation in \eqref{dualminmax} can be alternatively represented by solving its inherent maximization problem.

It's worth noting that in this maximization problem, $U(\lambda)$ adopts a hyperbox-metric form, as described in \eqref{hyperbox}. Therefore, $\xi$ can be further expressed as $\xi=\bar{d} + \lambda \Delta d z^+ -\lambda \Delta d z^-$, with both $z^+$ and $z^-$ being binary vectors indicating deviation direction. As outlined in \cite{fan2021flexibility}, the optimal $\xi$ must be achieved at the boundary of $U(\lambda)$.

For notation brevity, we suppose that $\xi\in\mathbb{R}^{N\times 1}$, $H_1\in\mathbb{R}^{M_1\times N}$, and $H_2\in\mathbb{R}^{M_2\times N}$. Given these, the expanded form of the maximization problem in \eqref{dualminmax} can be reformulated as:
\begin{subequations}\label{RDFEA}
\begin{align}
\psi &=\max_{z, \mu, \nu} \ \  h_1^T\mu + h_2^T\nu + \sum_{n=1}^N \sum_{m=1}^{M_1} (\bar{d}_n H_{1,m,n} \mu_m  \nonumber  \\
&+ \lambda \Delta d_n H_{1,m,n} \hat{\mu}^+_{n,m}-\lambda \Delta d_n H_{1,m,n} \hat{\mu}^-_{n,m}) \nonumber \\
&+ \sum_{n=1}^N \sum_{m=1}^{M_2} \big[\bar{d}_n H_{2,m,n} (\nu^a_m-\nu^b_m)  \nonumber  \\
&+ \lambda \Delta d_n H_{2,m,n} (\hat{\nu}^{a,+}_{n,m}-\hat{\nu}^{b,+}_{n,m}) \nonumber \\
&-\lambda\Delta d_n H_{2,m,n} (\hat{\nu}^{a,-}_{n,m}-\hat{\nu}^{b,-}_{n,m})\big] \nonumber\\
&-\sum_{n=1}^N(\bar{d}_n \gamma_n + \lambda \Delta d_n \gamma_n z^{+}_n - \lambda \Delta d_n \gamma_n z^{-}_n) \label{sub_obj}\\
\text{s.t.} \ \ & \text{Constraints}\ \eqref{mu_c1}-\eqref{mu_c2}, \\
&-z^+_n\le \hat{\mu}^+_{n,m}, \ \mu_m \le \hat{\mu}^+_{n,m} \le 1-z^+_n+\mu_m, \nonumber\\
&-z^-_n\le \hat{\mu}^-_{n,m}, \ \mu_m \le \hat{\mu}^-_{n,m} \le 1-z^-_n+\mu_m, \nonumber\\
& -1\le \hat{\mu}^+_{n,m} \le 0, \ -1\le \hat{\mu}^-_{n,m} \le 0, \nonumber\\ 
&\forall n=1\dots N, \ \forall m=1\dots M_1.\label{linearcons}\\
&-z^+_n\le \hat{\nu}^{\kappa,+}_{n,m}, \ \nu_m^\kappa \le \hat{\nu}^{\kappa,+}_{n,m} \le 1-z^+_n+\nu_m^\kappa, \nonumber\\
&-z^-_n\le \hat{\nu}^{\kappa,-}_{n,m}, \ \nu_m^\kappa \le \hat{\nu}^{\kappa,-}_{n,m} \le 1-z^-_n+\nu_m^\kappa, \nonumber\\
& -1\le \hat{\nu}^{\kappa,+}_{n,m}\le 0, \ -1\le \hat{\nu}^{\kappa,-}_{n,m}\le 0, \nonumber\\
&\forall \kappa \in \{a,b\}, \ \forall n=1\dots N, \ \forall m=1 \dots M_2.\\
&z^+_n+z^-_n= 1,\ \ \forall z^+_n,z^-_n\in\{0,1\}.
\end{align}    
\end{subequations}

To tackle the bilinear term $\xi^T H_1^T \mu$ in the objective function, we introduce auxiliary variables $\hat{\mu}^+_{n,m}$ and $\hat{\mu}^-_{n,m}$ to denote the products $z^+_n \mu_m$ and $z^-_n \mu_m$, respectively. For the term $\xi^T H_2^T \nu$, we decompose $\nu$ into $\nu^a-\nu^b$. Both these components, $\nu^a$ and $\nu^b$, are restricted to the range $[-1,0]$. We then apply a method similar to \eqref{linearcons} to linearize the expressions $\xi^T H_2^T \nu^a$ and $\xi^T H_2^T \nu^b$. Hence, the maximization part in \eqref{dualminmax} is reformulated into a mixed-integer linear programming model as \eqref{RDFEA}.

Upon solving \eqref{RDFEA}, the optimal solutions are denoted as $(z^*, \mu^*, \nu^*)$ with the corresponding optimal value of $\psi^*$. In accordance with \eqref{dualminmax}, we examine whether the following condition is satisfied:

\begin{equation}\label{checkcons}
\bar{d}^{T}\gamma \ + \psi^* \le \beta.
\end{equation}
If \eqref{checkcons} is met, the optimal solution to the master problem \eqref{DRO_OF}, denoted as $\lambda^*$, becomes the final flexibility result.

Otherwise, we refine the master problem by incorporating a feasibility cut $\bar{d}^{T}\gamma + \psi(\lambda,\gamma)\le \beta$. Here, $\psi(\lambda,\gamma)$ is derived by replacing with the optimal solution $(z^*, \mu^*, \nu^*)$ from \eqref{sub_obj}. Subsequently, the reformed master problem is developed as:
\begin{subequations}\label{finalmas}
\begin{align}
\max_{0\le\lambda\le 1} \ \ & \lambda \label{convexmas}\\
\text{s.t.} \ \ &  \sum_{n=1}^N \bigg\{\sum_{m=1}^{M_1} (\Delta d_n H_{1,m,n} \hat{\mu}^{+,*}_{n,m}-\Delta d_n H_{1,m,n} \hat{\mu}^{-,*}_{n,m})  \nonumber\\
&+ \sum_{m=1}^{M_2} \big[ \Delta d_n H_{2,m,n} (\hat{\nu}^{a+,*}_{n,m}-\hat{\nu}^{b+,*}_{n,m}) \nonumber\\
&- \Delta d_n H_{2,m,n} (\hat{\nu}^{a-,*}_{n,m}-\hat{\nu}^{b-,*}_{n,m})\big]  \nonumber\\
&- (\Delta d_n  z^{+,*}_n-\Delta d_n z^{-,*}_n)\gamma_n\bigg\}\lambda  \nonumber\\
&+h_1^T\mu^*+ h_2^T\nu^*+ \sum_{n=1}^N \bigg\{\sum_{m=1}^{M_1} \bar{d}_n H_{1,m,n} \mu_m^* \nonumber \\
&+ \sum_{m=1}^{M_2} \bar{d}_n H_{2,m,n} (\nu^{a*}_m-\nu^{b*}_m)\bigg\}   \le \beta, \label{feasibility_cut}
\end{align}
\end{subequations}
where the bilinear term $\lambda\gamma_n$ from the feasibility cut \eqref{feasibility_cut} is substituted with $w_n$, as depicted in \eqref{MC_EQ1}-\eqref{MC_EQ2} using McCormick Envelopes. 
\begin{subequations}\label{McCormick}
\begin{align}
& w_n \ge -\lambda K,\ w_n \ge \gamma_n +\lambda K - K,\ \forall n=1\dots N, \label{MC_EQ1}\\
& w_n \le \gamma_n -\lambda K + K,\ w_n \le \lambda K,\ \forall n=1\dots N. \label{MC_EQ2}
\end{align}
\end{subequations}
Notably, $K$ is a sufficiently large constant, and the constraint $-K\le \gamma\le K$ provides relaxation for the unrestricted $\gamma$.

\vspace{0.2cm}
\noindent {\bf Solution Approach:} The subsequent steps outline the cutting plane algorithm used to resolve the DRO model \eqref{DRO_model}:
\begin{enumerate}[label={\arabic*.}]
    \item Solve the master problem \eqref{DRO_OF} to determine the optimal value, denoted as $\lambda^*$. \label{step1}
    \item Assess the feasibility of the subproblem by solving \eqref{RDFEA} with the obtained $\lambda^*$. 
    \item Evaluate the validity of condition \eqref{checkcons}:
        \begin{itemize}
            \item If it holds, conclude the process and yield both the optimal solution and the master problem's objective value, $\lambda^*$.
            \item If not, refine the master problem by incorporating the feasibility cut from \eqref{feasibility_cut} and revert to Step \ref{step1}
        \end{itemize}
\end{enumerate}

\section{Real-Time Economic Dispatch}\label{sec:ED}
In this section, we present the mathematical framework for a Real-Time Economic Dispatch (RTED) problem, accounting for flexible resources such as power generators and Energy Storage Systems (ESS). The model is designed to optimize the generation output ($p^{\text G}_{n,t}$) and the net power output from ESSs ($p^{\text{ESS}}_{i,t}$) for each time interval $t$. Notably, a negative $p^{\text{ESS}}_{i,t}$ indicates ESS charging, while a positive value signals discharging. With a predefined uncertainty space $\mathcal{U}$, the feasible domain for these decision variables can be expressed as:

\begin{subequations}\label{linear}
	\begin{align}
	& \bold{X}(d) = \Bigg\{ \text{P}_{n}^{\text{min}} \leq p^{\text G}_{n,t} \leq \text{P}_{n}^{\text{max}},\ \forall n \in \mathcal{{G}},  \label{eqn-ED:PGLBD}\\
	& -{RD}_n \le p^{\text G}_{n,t}-p^{\text G}_{n,t-1} \le {RU}_n, \ \forall n \in \mathcal{{G}},  \label{eqn-ED:RAMPBD} \\
        & E^{\text{ESS}}_{i,t} = E^{\text{ESS}}_{i,t-1} - p^{\text{ESS}}_{i,t}, \ \forall i \in \mathcal{E}, \label{eqn-ED:SOCchange} \\
        & \text{E}^{\text{min}}_{i} \le E^{\text{ESS}}_{i,t} \le \text{E}^{\text{max}}_{i}, \ \forall i \in \mathcal{E}, \label{eqn-ED:SOClimit} \\
        & -\text{P}^{\text{max}}_{c,i} \le p^{\text{ESS}}_{i,t} \le \text{P}^{\text{max}}_{dc,i}, \ \forall i \in \mathcal{E},  \label{eqn-ED:Plimit} \\
	& \bigg|\sum_{b \in \mathcal{B}} \mbox{SF}_{b,l}(\sum_{n \in \mathcal{G}^{b}}p^{\text G}_{n,t} +\sum_{i \in \mathcal{E}^{b}}p^{\text{ESS}}_{i,t} - {d}_{b,t})\bigg|\leq F_{l}, \ \forall l \in \mathcal{L}, \label{transmissioncap}\\
	& \sum_{t \in \mathcal{T}}\bigg\{\sum_{n \in \mathcal{{G}}} C^{\text G}_{n}p^{\text G}_{n,t} + \sum_{i \in \mathcal{E}}C^{\text {ESS}}_{i}p^{\text{ESS}}_{i,t}\bigg\} \leq \tau,  \label{eqn-ED:costfun} \\
	& \sum_{n\in {\mathcal{G}}}p^{\text G}_{n,t} + \sum_{i \in \mathcal{E}}p^{\text{ESS}}_{i,t} - \sum_{b \in \mathcal{B}}{d}_{b,t}=0, \ \forall d_{b,t}\in \mathcal{U}_{b,t}, \label{balance} \\
	&  p^{\text G}_{n,t} \geq 0, \ \forall n \in {\mathcal{G}}, \ \forall t \in \mathcal{T}\Bigg\}. \label{range}
	\end{align}
\end{subequations}

We denote the sets of generators, ESSs, and transmission lines as $\mathcal{G}$, $\mathcal{E}$, and $\mathcal{L}$, respectively, with $\mathcal{G}^b$ and $\mathcal{E}^b$ indicating subsets of generators and ESSs at bus $b$. This model bypasses simultaneous ESS charging and discharging scenarios for arbitrage, given that its absence doesn't compromise system flexibility.


\begin{figure}[!h]
\centering
\includegraphics[width=\linewidth]{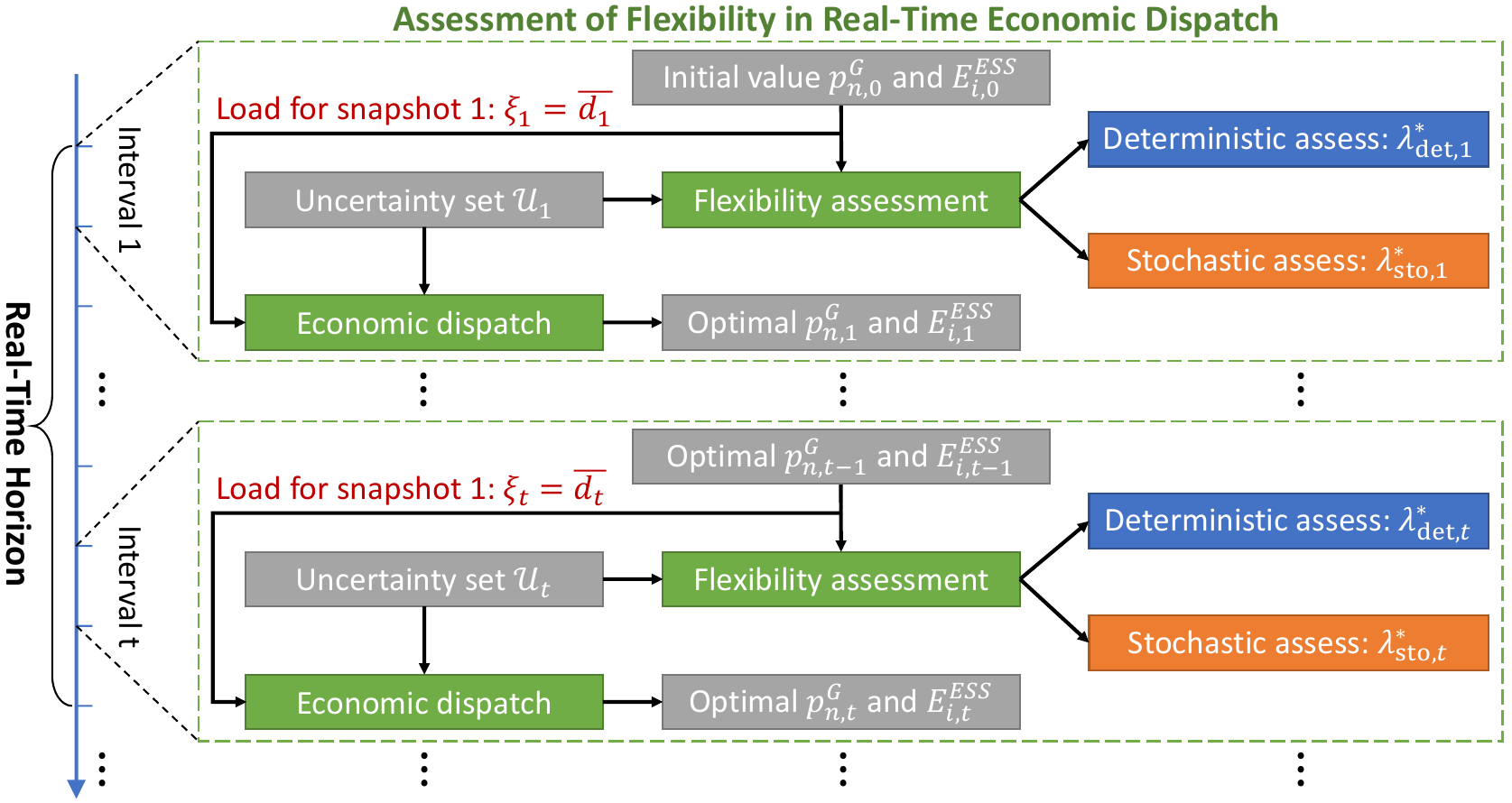}
\caption{Illustration of the flexibility assessment process for the RTED model.}
\label{fig1}
\end{figure}  

In practice, RTED is executed sequentially every 5–15 min interval with the static snapshot forecast data \cite{reddy2015real}. As depicted in Fig. \ref{fig1}, when optimizing for a given interval $t$, the decision variables from the previous interval $p^{\text G}_{n,t-1}$ and $E^{\text{ESS}}_{i,t-1}$, are treated as fixed constants, derived from the prior interval's economic dispatch results. With the net load's uncertainty set $\mathcal{U}_{t}$, we then optimize the flexibility parameter $\lambda_t$ for each interval $t$ using either deterministic or stochastic assessment, subject to the interval-specific constraints outlined in \eqref{linear}.

\section{Numerical Experiments}
In this section, we evaluate the maximum extent of net load uncertainty across all buses within the RTED model presented in Section \ref{sec:ED}. We employ both deterministic and stochastic assessments, applying them to four IEEE standard systems. This sequential RTED process is conducted at 5-minute intervals, over a defined 120-minute scheduling window.

\subsection{Flexibility Metric}
\begin{figure}[!h]
\centering
\includegraphics[width=\linewidth]{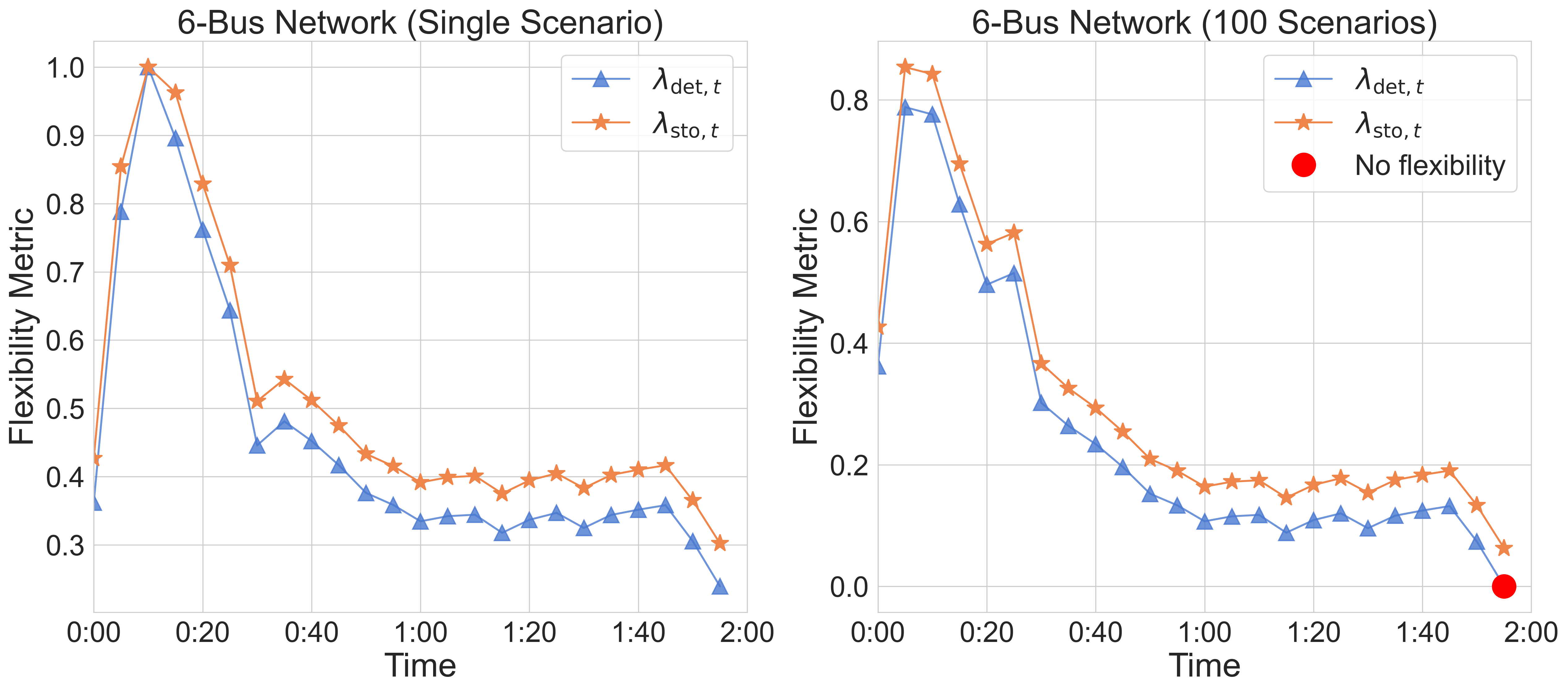}
\caption{Comparative flexibility outcomes: deterministic vs. stochastic assessments across single and multi net-load scenarios.}
\label{fig2}
\end{figure}

Fig. \ref{fig2} displays the RTED flexibility outcomes through deterministic and stochastic assessments. We present results from both single-scenario assessments, based on the nominal net loads, and those derived from 100-scenario assessments. The majority of the 100 scenarios are distributed within a narrow [0.99, 1.01] range relative to the normal net load for each time interval. Nevertheless, we incorporated an outlier scenario at 1.09 times the normal level to examine the responses of both deterministic and stochastic assessments to rare extreme cases in the power system. Notably, $\lambda_{\text{sto},t}$ consistently outperforms $\lambda_{\text{det},t}$ in both subplots. This difference arises from the DRO model's allowance in the stochastic assessment to accommodate minor system constraint violations. Specifically, we set $\beta$ in \eqref{DRO_CONS} to $0.05$, signifying the expected system operation violation below 5\%, thus enhancing the system's adaptability. Moreover, as the number of scenarios expands, both the deterministic flexibility $\lambda_{\text{det},t}$ and the stochastic flexibility $\lambda_{\text{sto},t}$ diminish in the right subplot. Influenced by extreme cases, $\lambda_{\text{det},t}$ decreases to zero in the final interval, suggesting the system lacks flexibility at that point. In contrast, $\lambda_{\text{sto},t}$ retains a flexibility measure of $0.063$, representing the flexibility exhibited in most scenarios, barring the extreme one.

\subsection{Sensitivity Analysis}
\begin{figure}[!h]
\centering
\includegraphics[width=\linewidth]{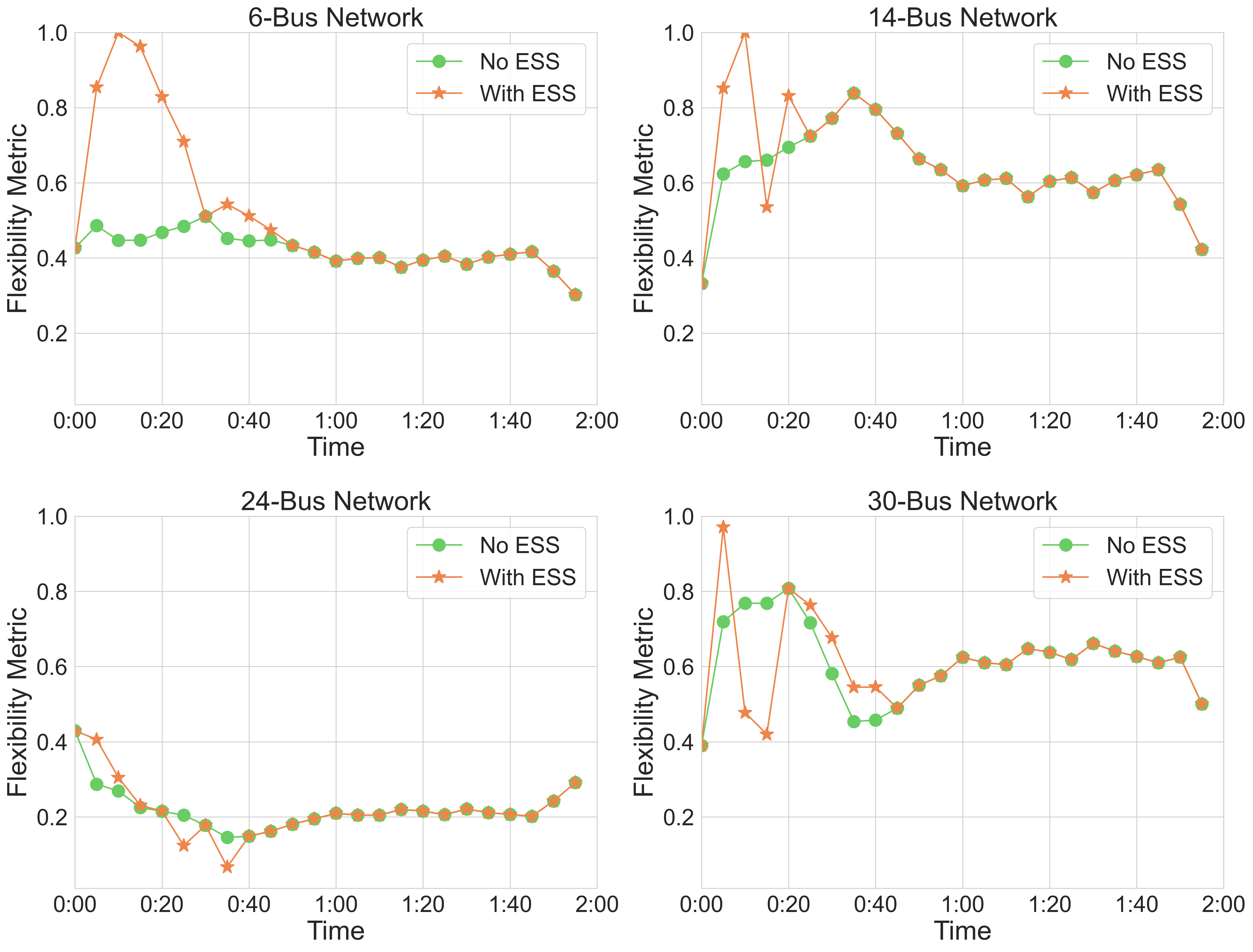}
\caption{ESS impact on system flexibility metric under stochastic assessment.}
\label{fig3}
\end{figure}

We investigate the impact on the system flexibility metric when there is a presence or absence of ESSs in the RTED model \eqref{linear}. By incorporating or excluding ESS-related constraints, specifically \eqref{eqn-ED:SOCchange}-\eqref{eqn-ED:Plimit}, and optimizing the system's real-time flexibility, we illustrate the comparative outcomes of $\lambda_{\text{sto},t}$ with and without ESSs in Fig. \ref{fig3}.

During the initial intervals of the scheduling window, power systems with ESSs demonstrate greater flexibility than those without. However, this advantage lessens over time, and sometimes, systems with ESSs can be less flexible. This is because as the ESS discharges, system generators curtail their power output in earlier intervals. Given the ramping rate constraint, power generation in later intervals may fall behind that of systems without ESSs, reducing flexibility when the ESS approaches the minimum charge. This trend is more evident in larger networks, i.e., the 24-bus and 30-bus systems, where the vast generation capacity diminishes the initial flexibility gains from ESSs.

\subsection{Computation Performance}
To verify the efficiency of the McCormick relaxation technique in solving the stochastic flexibility assessment $\lambda_{\text{sto},t}$, we also address the master problem \eqref{finalmas} with the nonconvex feasibility cut directly in Gurobi, referred to as the Gurobi-NC method. We limit the stochastic cutting plane algorithm to 30 iterations for each time interval due to time constraints. Table \ref{tab1} displays the computational times and the convergence performance for both approaches when handling the RTED with ESSs. Specifically, the convergence metric measures the number of intervals that a given method converges within 30 iterations, out of a total of 24 within the scheduling window.
\begin{table}[htbp]
  \centering
    \aboverulesep=0ex 
    \belowrulesep=0ex 
  \caption{Running time and convergence performance for the Gurobi-NC and McCormick method under stochastic assessment}
    \begin{tabular}{lcc|cc}
    \toprule
    \multicolumn{1}{c}{\multirow{2}[4]{*}{Networks}} & \multicolumn{2}{c|}{Time (seconds)} & \multicolumn{2}{c}{Convergence Metric (intervals)} \\
\cmidrule{2-5}          & \multicolumn{1}{l}{Gurobi-NC} & \multicolumn{1}{l|}{McCormick} & \multicolumn{1}{l}{Gurobi-NC} & \multicolumn{1}{l}{McCormick} \\
    \midrule
    6-bus & 4.92  & 4.54  & 24  & 24 \\
    14-bus & 28.72 & 18.31 & 23  & 24 \\
    24-bus & 6853.59 & 1749.41 & 1  & 24 \\
    30-bus & 19601.20 & 3007.45 & 8  & 24 \\
    \bottomrule
    \end{tabular}%
  \label{tab1}%
\end{table}%

As the system size increases, the efficiency of the McCormick method in optimizing the flexibility metric surpasses the Gurobi-NC method. While both methods converge optimally for the 6-bus system under the preset stopping criteria, the Gurobi-NC method struggles to do so for larger test systems during certain time intervals. In contrast, the McCormick method consistently achieves convergence. Owing to its stable convergence and shorter computational time, the McCormick method stands out as the preferred choice for solving the stochastic flexibility assessment in the DRO model.


\section{Conclusion}
Building on the deterministic assessment of power system flexibility, this paper introduced a stochastic assessment framework within the DRO model, which was tested through an RTED problem. Numerical results indicated that our stochastic assessment yielded less conservative flexibility metrics. Through sensitivity analysis, we observed that with larger system scales, the presence of numerous generators reduced the positive impact of ESSs on system flexibility. Additionally, the efficiency of the McCormick envelope in solving the DRO model was confirmed against the direct nonconvex approach.

\bibliographystyle{IEEEtran}
\bibliography{scholar}

\vfill

\end{document}